\documentclass[12pt]{article}

\usepackage{amsfonts,amsmath,amssymb,bbm,graphicx,dsfont}
\usepackage{setspace} 

\def\proof{\noindent{\bf Proof:}\hskip10pt}        
\def\QED{\hfill $\Box$}

\font\tenmath=msbm10 scaled 1200
\font\sevenmath=msbm7 scaled 1200
\font\Fivemath=msbm5 scaled 1200
\newfam\mathfam \textfont\mathfam=\tenmath
\scriptfont\mathfam=\sevenmath \scriptscriptfont\mathfam=\Fivemath

\def \\ { \cr }
\def\R{\mathbb{R}}
\def \1{1 \mkern -6mu 1} 
\def\N{\mathbb{N}}
\def\E{\mathbb{E}}
\def\P{\mathbb{P}}
 
\def\p{\mathcal{P}}

\def\R{\mathbb{R}}

\def\T{{\mathbb{T}}}

\def \e{{\rm e}}

\def \X{{\mathbf X}}

\def \d{{\rm d}}

\def \dt{{\ensuremath{\textup{d}}}}

\def \pY{{\ensuremath{\mathcal{Y}}}}
\def \pZ{{\ensuremath{\mathcal{Z}}}}

\def \Fr{{\ensuremath{\textup{Fr}}}}
\def \G{{\ensuremath{\textup{G}}}}

\newtheorem{theorem}{Theorem} 
\newtheorem{proposition}{Proposition} \newtheorem{lemma}{Lemma}
\newtheorem{corollary}{Corollary} 

\begin{document}

\title{The fragmentation process of an infinite recursive tree
  \\
  and Ornstein-Uhlenbeck type processes} \author{{Erich
    Baur\footnote{erich.baur@math.uzh.ch} { and } Jean
    Bertoin\footnote{jean.bertoin@math.uzh.ch}}\\ ENS Lyon and
  Universit\"at Z\"urich}
\maketitle 
\thispagestyle{empty}
\begin{abstract}
  We consider a natural destruction process of an infinite recursive tree
  by removing each edge after an independent exponential time. The
  destruction up to time $t$ is encoded by a partition $\Pi(t)$ of $\N$
  into blocks of connected vertices. Despite the lack of exchangeability,
  just like for an exchangeable fragmentation process, the process $\Pi$ is
  Markovian with transitions determined by a splitting rates measure ${\bf
    r}$. However, somewhat surprisingly, ${\bf r}$ fails to fulfill the
  usual integrability condition for the dislocation measure of exchangeable
  fragmentations.  We further observe that a time-dependent normalization
  enables us to define the weights of the blocks of $\Pi(t)$.  We study the
  process of these weights and point at connections with Ornstein-Uhlenbeck
  type processes.
  
\end{abstract}

\section{Introduction}
The purpose of this work is to investigate various aspects of a simple and
natural fragmentation process on an infinite tree, which turns out to
exhibit nonetheless some rather unexpected features.

Specifically, we first construct a tree $\T$ with set of vertices
$\N=\{1,\ldots\}$ by incorporating vertices one after the other and
uniformly at random. That is, $1$ is the root, and for each vertex $i\geq
2$, we pick its parent $u_i$ according to the uniform distribution in
$\{1,\ldots, i-1\}$, independently of the other vertices. We call $\T$ an
infinite (random) recursive tree. Recursive trees are especially useful in
computer science where they arise as data structures; see e.g. the survey
by Mahmoud and Smythe \cite{MaSmy} for background.

We next destroy $\T$ progressively by removing each edge $e_i$ connecting
$i$ to its parent $u_i$ at time $\epsilon_i$, where the sequence
$(\epsilon_i: i\geq 2)$ consists of i.i.d. standard exponential variables,
which are further independent of $\T$. Panholzer \cite{Pan} investigated
costs related to this destruction process, whereas in a different
direction, Goldschmidt and Martin~\cite{GM} used it to provide a remarkable
construction of the Bolthausen-Sznitman coalescent. We also refer to Kuba
and Panholzer \cite{KP1} for the study of a related algorithm for isolation
of nodes, and to our survey \cite{BaBe} for further applications and many
more references.

 Roughly
speaking, we are interested here in the fragmentation process that results
from the destruction.  We represent the destruction of $\T$ up to time $t$
by a partition $\Pi(t)$ of $\N$ into blocks of connected vertices. In other
words, if we view the fragmentation of $\T$ up to time $t$ as a Bernoulli
bond-percolation with parameter $\e^{-t}$, then the blocks of $\Pi(t)$ are
the percolation clusters. Clearly $\Pi(t)$ gets finer as $t$ increases, and
it is easily seen from the fundamental splitting property of random
recursive trees that the process $\Pi= (\Pi(t): t\geq 0)$ is Markovian.  In
this direction, we also recall that Aldous and Pitman~\cite{AP} have
considered a similar logging of the Continuum Random Tree that yields a
notable fragmentation process, dual to the standard additive coalescent. We
further point at the very recent work~\cite{KBN} in which the effects of
repeated random removal of nodes (instead of edges) in a finite random
recursive tree are analyzed.

It turns out that $\Pi$ shares many features similar to homogeneous
fragmentation processes as defined in~\cite{Be1, RFCP}. In particular, the
transition kernels of $\Pi$ are very similar to those of a homogeneous
fragmentation; they are entirely determined by the so-called splitting
rates ${\bf r}$, which define an infinite measure on the space of
partitions of $\N$. However, there are also major differences:
exchangeability, which is a key requirement for homogeneous fragmentation
processes, fails for $\Pi$, and perhaps more notably, the splitting rates
measure ${\bf r}$ does not fulfill the fundamental integral
condition~\eqref{eqIC} which the splitting rates of homogeneous
fragmentation processes have to satisfy.

It is known from the work of Kingman~\cite{Ki} that exchangeability plays a
fundamental role in the study of random partitions, and more precisely, it
lies at the heart of the connection between exchangeable random partitions
(which are discrete random variables), and random mass-partitions (which
are continuous random variables). In particular, the distribution of an
exchangeable random partition is determined by the law of the asymptotic
frequencies of its blocks $B$,
\begin{equation}\label{eqasymf}
|B|=\lim_{n\to\infty} n^{-1}\#\{i\leq n : i\in B\}.
\end{equation}

Even though $\Pi(t)$ is not exchangeable for $t>0$, it is elementary to see
that every block of $\Pi(t)$, say $B(t)$, has an asymptotic frequency. However
this asymptotic frequency is degenerate, $|B(t)|=0$ (note that if $\Pi$
were exchangeable, this would imply that all the blocks of $\Pi(t)$ would
be singletons). We shall obtain a finer result and show that the limit
\begin{equation}\label{eqmass}
\lim_{n\to\infty} n^{-\e^{-t}}\#\{i\leq n : i\in B(t)\}
\end{equation}
exists in $(0,\infty)$ almost surely. 
We shall refer to the latter as the weight of the block $B(t)$ (we stress
that this definition depends on the time $t$ at which the block is taken),
and another natural question about the destruction of $\T$ is thus to
describe the process $\X$ of the weights of the blocks of the
partition-valued process $\Pi$.

Because $\Pi$ resembles homogeneous fragmentations, but with splitting rate
measure ${\bf r}$ which does not fulfill the integral condition of the
former, and because the notion~\eqref{eqmass} of the weight of a block
depends on the time $t$, one might expect that $\X$ should be an example of
a so-called compensated fragmentation which was recently introduced
in~\cite{Be2}. Although this is not exactly the case, we shall see that
$\X$ fulfills closely related properties.  Using well-known connections
between random recursive trees, Yule processes, and P\'olya urns,
cf.~\cite{BaBe}, we shall derive a number of explicit results about its
distribution. In particular, we shall show that upon a logarithmic
transform, $\X$ can be viewed as a branching Ornstein-Uhlenbeck process.

The rest of this paper is organized as follows. In Section 2, we study the
structure of the partition-valued process $\Pi$ which stems from the
destruction of $\T$, stressing the resemblances and the differences with
exchangeable fragmentations.  In Section 3, we observe that after a
suitable renormalization that depends on $t$, the blocks of the partition
$\Pi(t)$ possess a weight, and we relate the process of these weights to
Ornstein-Uhlenbeck type processes.

\section{Destruction of $\T$ and fragmentation of partitions}
The purpose of this section is to show that, despite the lack of
exchangeability, the partition valued process $\Pi$ induced by the
fragmentation of $\T$ can be analyzed much in the same way as a homogeneous
fragmentation. We shall present the main features and merely sketch proofs,
referring to Section 3.1 in~\cite{RFCP} for details.

We start by recalling that a partition $\pi$ of $\N$ is a sequence $(\pi_i:
i\in \N)$ of pairwise disjoint blocks, indexed in the increasing order of
their smallest elements, and such that $\sqcup_{i\in \N} \pi_i=\N$. We
write $\p$ for the space of partitions of $\N$, which is a compact
hypermetric space when endowed with the distance
$${\rm d}(\pi, \pi')= 1/\max\{n\in\N: \pi_{\mid [n]}= \pi'_{\mid [n]}\},$$
where $\pi_{\mid B}$ denotes the restriction of $\pi$ to a subset
$B\subseteq \N$ and $[n]=\{1, \ldots, n\}$ is the set of the $n$ first
integers. The space $\p_B$ of partitions of $B$ is defined similarly.

We next introduce some spaces of functions on $\p$. First, for every $n\geq
1$, we write $D_n$ for the space of functions $f: \p\to \R$ which remain
constant on balls with radius $1/n$, that is such that $f(\pi)=f(\eta)$
whenever the restrictions $\pi_{\mid [n]}$ and $\eta_{\mid [n]}$ of the
partitions $\pi$ and $\eta$ to $[n]$ coincide. Plainly, $D_n\subset
D_{n+1}$, and we set $$D_{\infty}=\bigcup_{n\geq 1} D_n.$$ Observe that
$D_{\infty}$ is a dense subset of the space ${\mathcal C}(\p)$ of
continuous functions on $\p$.

In order to describe a family of transition kernels which appear naturally
in this study, we first need some notation. For every block $B\subseteq\N$,
write $B(j)$ for the $j$-th smallest element of $B$ (whenever it makes
sense), and then, for every partition $\pi\in\p$, $B \circ \pi$ for the
partition of $B$ generated by the blocks $B(\pi_i)= \{B(j): j\in \pi_i\}$
for $i\in\N$. In other words, $B \circ \pi$ is simply the partition of $B$
induced by $\pi$ when one enumerates the elements of $B$ in their natural
order. Of course, if the cardinality of $B$ is finite, say equal to $k\in\N$,
then $B \circ \pi$ does only depend on $\pi$ through $\pi_{\mid[k]}$, so
that we may consider $B \circ \pi$ also for $\pi\in\p_{[k]}$ (or $\pi\in
\p_{[\ell]}$ for any $\ell\geq k$).

In the same vein, for partitions $\eta\in\p_B$ and every integer $i\geq 1$,
we write $\eta{\underset{i}{\,\circ\,}}\pi$ for the partition of $B$ that
results from fragmenting the $i$-th block of $\eta$ by $\pi$, that is
replacing the block $\eta_i$ in $\eta$ by $\eta_i \circ \pi$. Again, if
$k=\#B<\infty$, we may also take $\pi\in\p_{[\ell]}$ for $\ell\geq k$.

Finally, for every $k\geq 2$, we consider a random partition of $\N$ that
arises from the following P\'olya urn.  At the initial time, the urn
contains $k-1$ black balls labeled $1, \ldots, k-1$ and a single red ball
labeled $k$. Balls with labels $k+1, k+2, \ldots$ are colored black or red
at random and then incorporated to the urn one after the other. More
precisely, for $n\geq k$, the color given to the $n+1$-th ball is that of a
ball picked uniformly at random when the urn contains $n$ balls. This
yields a random binary partition of $\N$; we write ${\bf p}_k$ for its
law. We set
\begin{equation}\label{defbfr}
{\bf r}=\sum_{k=2}^{\infty} {\bf p}_k,
\end{equation}
which is thus an infinite measure on the set of binary partitions of $\N$. 

Recall that each edge of $\T$ is deleted at an exponentially distributed
random time, independently of the other edges. This induces, for every
$t\geq 0$, a random partition $\Pi(t)$ of $\N$ into blocks corresponding to
the subsets of vertices which are still connected at time $t$. Observe
that, by construction and the very definition of the distance on $\p$, the
process $\Pi$ has c\`adl\`ag paths.

We are now able to state the main result of this section.

\begin{theorem}\label{T1} {\rm (i)} The process $\Pi=(\Pi(t): t\geq 0)$ is
  Markovian and has the Feller property. We write $\G$ for its infinitesimal
  generator.

  \noindent{\rm (ii)} For every $n\geq 1$, $D_n$ is invariant and therefore
  $D_{\infty}$ is a core for $\G$.

  \noindent{\rm (iii)} For every $f\in D_{\infty}$ and $\eta\in\p$, we have
$$\G f(\eta) = \int_{\pi\in\p}{\bf r}(\d \pi) \sum_i \left(f(\eta{\underset{i}{\,\circ\,}}\pi)-f(\eta)\right).$$
\end{theorem}

We stress that this characterization of the law of the process $\Pi$ is
very close to that of a homogeneous fragmentation. Indeed, one can rephrase
well-known results (cf. Section 3.1.2 in~\cite{RFCP}) on the latter as
follows. Every homogeneous fragmentation process $\Gamma=(\Gamma_t: t\geq
0)$ is a Feller process on $\p$, such that the sub-spaces $D_n$ are
invariant (and hence $D_{\infty}$ is a core). Further, its infinitesimal
generator A is given in the form
$$\textup{A}f(\eta) = \int_{\pi\in\p}{\bf s}(\d \pi) \sum_i \left(f(\eta{\underset{i}{\,\circ\,}}\pi)-f(\eta)\right)$$
for every $f\in D_{\infty}$ and $\eta\in \p$, where ${\bf s}$ is some {\it
  exchangeable} measure on $\p$. More precisely, ${\bf s}(\{{\bf
  1}_{\N}\})=0$, where for every block $B\subseteq \N$, ${\bf 1}_{B}\in
\p_B$ denotes the neutral partition which has a single non-empty block $B$,
and
\begin{equation}\label{eqbfs}
{\bf s}\left(\left\{\pi\in\p: \pi_{\mid [n]}\neq {\bf 1}_{[n]}\right\}\right)<\infty\qquad \hbox{for all }n\geq 2.
\end{equation}
Observe that the measure ${\bf r}$ fails to be
exchangeable, but it fulfills~\eqref{eqbfs}; indeed, one has
$$
{\bf r}\left(\left\{\pi\in\p: \pi_{\mid [n]}\neq {\bf
      1}_{[n]}\right\}\right) = \sum_{k=2}^{n} {\bf p}_k(\p) = n-1.
$$

We shall now prepare the proof of Theorem~\ref{T1}. In this direction, it
is convenient to introduce some further notation.  Consider an arbitrary
block $B\subseteq \N$, a partition $\eta\in \p_B$ and a sequence
$\pi^{(\cdot)}=(\pi^{(i)}: i\in\N)$ in $\p$. We write $ \eta \circ
\pi^{(\cdot)} $ for the partition of $B$ whose family of blocks is given by
those of $\eta {\underset{i}{\,\circ\,}}\pi^{(i)}$ for $i\in\N$.  In words,
for each $i\in\N$, the $i$-th block of $\eta$ is split according to the
partition $\pi^{(i)}$.  Next, consider a probability measure ${\bf q}$ on
$\p$ and a sequence $(\pi^{(i)}: i\in\N)$ of i.i.d. random partitions with
common law ${\bf q}$. We associate to ${\bf q}$ a probability kernel
$\Fr(\cdot, {\bf q})$ on $\p_B$, by denoting the distribution of $ \eta
\circ \pi^{(\cdot)} $ by $\Fr(\eta, {\bf q})$ for every $\eta \in\p_B$.  We
point out that if ${\bf q}$ is exchangeable, then $ \eta \circ
\pi^{(\cdot)} $ has the same distribution as the random partition whose
blocks are given by the restrictions $\pi^{(i)}_{\mid \eta_i}$ of
$\pi^{(i)}$ to $ \eta_i$ for $i\in\N$, and $\Fr(\cdot, {\bf q})$ thus
coincides with the fragmentation kernel that occurs for homogeneous
fragmentations, see Definition 3.2 on page 119 in \cite{RFCP}. Of course,
the assumption of exchangeability is crucial for this identification to
hold.

Note that the restriction of partitions of $\N$ to $[n]$ is compatible with the
fragmentation operator $\Fr(\cdot,\cdot)$, in the sense that
\begin{equation}
\label{eqcompatible}
\left(\eta\circ\pi^{(\cdot)}\right)_{\mid[n]} = \eta_{\mid[n]}\circ{\pi^{(\cdot)}}_{\mid[n]}.
\end{equation}

\begin{proposition}\label{P1} For every $n\in \N$, the
  process $\Pi_{\mid [n]}=(\Pi_{\mid [n]}(t): t\geq 0)$ obtained by
  restricting $\Pi$ to $[n]$, is a continuous time Markov chain on
  $\p_{[n]}$.
 
  Its semigroup can be described as follows: for every $s,t\geq0$, the
  conditional distribution of $\Pi_{\mid [n]}(s+t)$ given $\Pi_{\mid
    [n]}(s)=\eta$ is $\Fr(\eta, {\bf q}_t)$, where ${\bf q}_t$ denotes the
  distribution of $\Pi(t)$.
\end{proposition}

\proof The case $n=1$ is clear since $\Pi_{\mid [1]}(t)=(\{1\},\emptyset,\dots)$ for all times
$t\geq 0$. Assume now $n\geq 2$. The proof
relies crucially on the so-called splitting property of random recursive
trees that we now recall (see, e.g.  Section 2.2 of~\cite{BaBe}). Given a
subset $B\subseteq \N$, the image of $\T$ by the map $j\mapsto B(j)$ which
enumerates the elements of $B$ in the increasing order, is called a random
recursive tree on $B$ and denoted by $\T_B$. In particular, for $B=[n]$,
the restriction of $\T$ to the first $n$ vertices is a random recursive
tree on $[n]$. Imagine now that we remove $k$ fixed edges (i.e. edges with
given indices, say $ i_1, \ldots, i_k$, where $2\leq i_1< \ldots< i_k\leq
n$) from $\T_{ [n]}$. Then, conditionally on the induced partition of
$[n]$, say $\eta=(\eta_1, \ldots, \eta_{k+1})$, the resulting $k+1$
subtrees are independent random recursive trees on their respective sets of
vertices $\eta_j$, $j=1,\ldots, k+1$.

It follows easily from the lack of memory of the exponential distribution
and the compatibility property~\eqref{eqcompatible} that the restricted
process $\Pi_{\mid [n]}= \left( \Pi_{\mid [n]}(t): t\geq 0\right)$ is a
continuous time Markov chain. More precisely, the conditional distribution
of $\Pi_{\mid [n]}(s+t)$ given $\Pi_{\mid [n]}(s)=\pi$ is $\Fr(\pi, {\bf
  q}_t)$, where $\Fr(\cdot,{\bf q}_t)$ is here viewed as a probability
kernel on $\p_{[n]}$.\QED
 
In order to describe the infinitesimal generator of the restricted
processes $\Pi_{\mid [n]}$ for $n\in\N$, we consider its rates of jumps,
which are defined by
$$r_{\pi} =\lim_{t\to 0+} t^{-1}\P(\Pi_{\mid [n]}(t)=\pi),$$
where now $\pi$ denotes a generic partition of $[n]$ which has at least two
(non-empty) blocks.  The rates of jumps $r_{\pi}$ determine the
infinitesimal generator $\G_n$ of the restricted chain $\Pi_{\mid [n]}$,
specifically we have for $f: \p_{[n]}\to \R$ and $\eta \in \p_{[n]}$
$$\G_nf(\eta) = \sum_{\pi\in\p_{[n]}}\sum_i \left(f(\eta{\underset{i}{\,\circ\,}}\pi)-f(\eta)\right) r_{\pi}$$
(recall that $\eta{\underset{i}{\,\circ\,}}\pi$ denotes the partition that
results from fragmenting the $i$-th block of $\eta$ according to $\pi$).
This determines the distribution of the restricted chain $\Pi_{\mid [n]}$,
and hence, letting $n$ vary in $\N$, also characterizes the law of
$\Pi$. Recall also that the measure ${\bf r}$ on $\p$ has been defined by
\eqref{defbfr}.

\begin{proposition} \label{P2} For every $n\geq 2$ and every partition $\pi$ of
  $[n]$ with at least two (non-empty) blocks, there is the identity
$$r_{\pi}={\bf r}(\p_{\pi}),$$
where $\p_{\pi}=\{\eta\in\p: \eta_{\mid [n]}=\pi\}$.
\end{proposition} 

\proof This should be intuitively straightforward from the connection
between the construction of random recursive trees and the dynamics of
P\'olya urns. Specifically, fix $n\geq 2$ and consider a partition
$\pi\in\p_{[n]}$.  If $\pi$ consists in three or more non-empty blocks,
then we clearly have
$$\lim_{t\to 0+} t^{-1}\P(\Pi_{\mid [n]}(t)=\pi)=0,$$
since at least two edges have to be removed from $\T_{\mid [n]}$ in order
to yield a partition with three or more blocks. Assume now that $\pi$ is
binary with non-empty blocks $\pi_1$ and $\pi_2$, and let $k=\min \pi_2$.
Then only the removal of the edge $e_k$ may possibly induce the partition
$\pi$, and more precisely, if we write $\eta$ for the random partition of
$[n]$ resulting from the removal of $e_k$, then the probability that
$\eta=\pi$ is precisely the probability that in a P\' olya urn containing
initially $k-1$ black balls labeled $1, \ldots, k-1$ and a single red ball
labeled $k$, after $n-k$ steps, the red balls are exactly those with
labels in $\pi_2$. Since the edge $e_k$ is removed at unit rate, this gives
$$\lim_{t\to 0+} t^{-1}\P(\Pi_{\mid [n]}(t)=\pi)= r_{\pi}={\bf p}_k(\p_{\pi})$$
in the notation of the statement.  Note that the right-hand side can be
also written as ${\bf r}(\p_{\pi})$, since ${\bf p}_\ell(\p_{\pi})=0$ for
all $\ell \neq k$.  \QED

Proposition~\ref{P2} should be compared with Proposition 3.2 in
\cite{RFCP}; we refer henceforth to ${\bf r}$ as the splitting rate of
$\Pi$.

We have now all the ingredients necessary to establish Theorem~\ref{T1}.

{\noindent{\bf Proof of Theorem~\ref{T1}:}\hskip10pt} From
Proposition~\ref{P1}, we see that the transition semigroup of $\Pi$ is
given by $\left(\Fr(\cdot, {\bf q}_{t}): t\geq 0\right)$, and it is easily
checked that the latter fulfills the Feller property; cf. Proposition
3.1(i) in~\cite{RFCP}. Point (ii) is immediate from the compatibility of
restriction with the fragmentation operator,
see~\eqref{eqcompatible}. Concerning (iii), let $f\in D_n$. Since $f$ is
constant on $\{\eta\in\p: \eta_{\mid [n]}=\pi\}$, it can naturally be
restricted to a function $f:\p_{[n]}\rightarrow\mathbb{R}$.  By the
compatibility property~\eqref{eqcompatible}, with $r_{\pi'}={\bf
  r}(\p_{\pi'})$ for a partition $\pi'$ of $[n]$ defined as in
Proposition~\ref{P2}, we obtain
$$\int_{\pi\in\p}{\bf r}(\d \pi) \sum_i
\left(f(\eta{\underset{i}{\,\circ\,}}\pi)-f(\eta)\right) =
\sum_{\pi'\in\p_{[n]}}\sum_i
\left(f(\eta_{\mid[n]}{\underset{i}{\,\circ\,}}\pi')-f(\eta_{\mid[n]})\right)
r_{\pi'}=\G_nf(\eta_{\mid[n]}),$$ where $\G_n$ is the infinitesimal generator
of the restricted chain $\Pi_{\mid[n]}$ found above. This readily yields (iii). \QED

\noindent{\bf Remark.} It may be interesting to recall that the standard
exponential law is invariant under the map $t\mapsto -\ln(1-\e^{-t})$, and
thus, if we set $\hat \epsilon_i = -\ln(1-\exp(-\epsilon_i))$ (recall that
$\epsilon_i$ is the instant at which the edge connecting the vertex $i$ to
its parent is removed), then $(\hat \epsilon_i)_{i\geq 2}$ is a sequence of
i.i.d. exponential variables. The time-reversal $t\mapsto -\ln(1-\e^{-t})$
transforms the destruction process of $\T$ into a construction process of
$\T$ defined as follows. At each time $\hat \epsilon_i $, we create an edge
between $i\geq 2$ and its parent which is chosen uniformly at random in
$\{1,\ldots, i-1\}$. It follows that the time-reversed process $\hat \Pi(t)
= \Pi(-\ln(1-\e^{-t})-)$, $t\geq 0$, is a binary coalescent process such
that the rate at which two blocks, say $B$ and $B'$ with $\min B < \min
B'$, merge, is given by $\#\{j\in B: j<\min B'\}$. This can be
viewed as a duality relation between fragmentation and coalescent
processes; see Dong, Goldschmidt and Martin~\cite{DGM} and references
therein.

Recall that for every
$k\geq 2$, a random binary partition with law ${\bf p}_k$ resulting from
the P\'olya urn construction possesses asymptotic frequencies in the
sense of~\eqref{eqasymf}. This readily entails the following result.

\begin{proposition}\label{P3} For ${\bf r}$-almost all binary partitions
  $(B_1,B_2)\in\p$, the blocks $B_1$ and $B_2$ have asymptotic frequencies,
  and more precisely, we have
$$\int_{\p}f(|B_1|,|B_2|)\d {\bf r} = \int_0^{1} f(1-x,x) x^{-2}\d x,$$
where $f:[0,1]^2\to \R_+$ denotes a generic measurable function. In particular,
$$ \int_{\p}(1-|B_1|)\d {\bf r} =\infty.$$
\end{proposition}
\proof Indeed, it is a well-known fact of P\'olya urns that for each $k\geq
2$, ${\bf p}_k$-almost every partition $(B_1,B_2)$ has asymptotic
frequencies with $|B_1|+|B_2|=1$ and $|B_2|$ has the beta distribution with
parameters $(1,k-1)$, i.e. with density $(k-1)(1-x)^{k-2}$ on $(0,1)$. Our
claims follow immediately since
$$\sum_{k=2}^{\infty} (k-1)(1-x)^{k-2} = x^{-2},\qquad x\in(0,1).$$\QED

It is interesting to recall that the splitting rates ${\bf s}$ of a
homogeneous fragmentation must fulfill the integrability condition
\begin{equation}\label{eqIC}
\int_{\p} (1-|B_1|) \d {\bf s} < \infty,
\end{equation}
which thus fails for ${\bf r}$ !

We next turn our attention to the Poissonian structure of the process
$\Pi$, which can be rephrased in terms similar to those in Section 3.1 of
\cite{RFCP}. In this direction, we introduce a random point measure
$${M}= \sum_{i=2}^{\infty} \delta_{(\epsilon_i, \Delta_i, k_i)}$$ on
$\R_+\times \p\times \N$ as follows. Recall that $\epsilon_i$ is the time
at which the edge $e_i$ connecting the vertex $i\in\N$ to its parent in
$\T$ is removed. Immediately before time $\epsilon_i$, the vertex $i$
belongs to some block of the partition $\Pi(\epsilon_i-)$, we denote the
label of this block by $k_i$ (recall that $\Pi$ is c\`adl\`ag and that
blocks of a partition are labeled in the increasing order of their smallest
element).  Removing the edge $e_i$ yields a partition of that block
$B=\Pi_{k_i}(\epsilon_i-)$ into two sub-blocks, which can be expressed
(uniquely) in the form $B\circ \Delta_i$. This defines the binary partition
$\Delta_i$ and hence the point measure ${M}$ unambiguously.  The process
$\Pi$ can be recovered from ${M}$, in a way similar to that explained on
pages 97-98 in~\cite{RFCP}. Roughly speaking, for every atom of $M$, say
$(t,\Delta,k)$, $\Pi(t)$ results from partitioning the $k$-th block of
$\Pi(t-)$ using $\Delta$, that is by replacing $\Pi_k(t-)$ by
$\Pi_k(t-)\circ \Delta$.  Adapting the arguments of Section 3.1.3
in~\cite{RFCP}, we have the following result:

\begin{proposition}\label{P4} The random measure ${M}$ is Poisson with intensity
  $\lambda \otimes {\bf r}\otimes \#$, where $\lambda$ denotes the Lebesgue
  measure on $\R_+$ and $\#$ the counting measure on $\N$.
\end{proposition}
\proof Recall that we write ${\bf 1}_{[n]}=([n],\emptyset,\ldots)$ for the
partition of $[n]$ which consists of a single non-empty block. Consider a
Poisson random measure $M'$ with intensity $\lambda \otimes {\bf r}\otimes
\#$ as in the statement.  Then $M'$ has almost surely at most one atom in
each fiber $\{t\}\otimes\p\otimes\N$, and the discussion
below~\eqref{eqbfs} shows that for each $t'\geq 0$ and every $n\in\N$, the
number of atoms $(t,\pi,k)$ of $M'$ with $t\leq t'$, $\pi_{\mid[n]}\neq
{\bf 1}_{[n]}$ and $k\leq n$ is finite. We may therefore define for fixed
$n\in\N$ a $\p_{[n]}$-valued continuous time Markov chain
$(\Pi'^{[n]}(t):t\geq 0)$ starting from $\Pi'^{[n]}(0) = {\bf 1}_{[n]}$ as
follows: If $t$ is a time at which the fiber $\{t\}\times\p\times\N$
carries an atom $(t,\pi,k)$ of $M'$ such that $\pi_{\mid[n]}\neq {\bf
  1}_{[n]}$ and $k\leq n$, then $\Pi'^{[n]}(t)$ results from
$\Pi'^{[n]}(t-)$ by replacing its $k$-th block $\Pi_k'^{[n]}(t-)$ by
$\Pi_k'^{[n]}(t-)\circ\pi_{\mid[n]}$.

The sequence $(\Pi'^{[n]}(t):n\in\N)$ is clearly compatible for every
$t\geq 0$, in the sense that ${\Pi'^{[n]}}_{\mid[m]}(t)=\Pi'^{[m]}(t)=$ for
all integers $n\geq m$. We deduce as in the proof of Lemma 3.3
in~\cite{RFCP} that there exists a unique $\p$-valued c\`adl\`ag function
$(\Pi'(t):t\geq 0)$ such that $\Pi'_{\mid [n]}(t) =
\Pi'^{[n]}(t)$. Moreover, the $i$-th block $\Pi'_i(t)$ of $\Pi'(t)$ is
given by the increasing union $\Pi'_i(t) =\cup_ {n\in\N}\Pi'^{[n]}_i(t)$,
and it follows from the very construction of $\Pi'^{[n]}(t)$ that the
process $\Pi'$ can be recovered from $M'$ similarly to the description above
the statement of the proposition.  It remains to check that $\Pi'$ and
$\Pi$ have the same law, which follows if we show that the restricted
processes $\Pi'_{\mid[n]}=\Pi'^{[n]}$ and $\Pi_{\mid [n]}$ have the same
law for each $n\in \N$. Fix $n\geq 2$, and denote by $\pi$ a partition of
$[n]$ with at least two non-empty blocks. From the Poissonian construction
of $\Pi'^{[n]}$, with $\p_{\pi}$ as in the statement of
Proposition~\ref{P2}, we first see that
$$
\lim_{t\to 0+} t^{-1}\P\left(\Pi'^{[n]}(t)=\pi\right)={\bf r}(\p_{\pi}).$$
Next, if $\pi'\neq \pi''\in\p_{[n]}$, the jump rate of $\Pi'^{[n]}$ from
$\pi'$ to $\pi''$ is non-zero only if $\pi''$ can be obtained from $\pi'$
by replacing one single block of $\pi'$, say the $k$-th block $\pi'_k$, by
$\pi'_k\circ\pi$, where $\pi$ is some binary partition of $[n]$. This
observation and the last display readily show that $\Pi'^{[n]}$ and
$\Pi_{\mid [n]}$ have the same generator, and hence their laws agree.\QED

\section{The process of the weights}
\label{SW}
Even though the splitting rates measure ${\bf r}$ of the fragmentation
process $\Pi$ fails to fulfill the integral condition~\eqref{eqIC}, we
shall see that we can nonetheless define the weights of its blocks.  The
purpose of this section is to investigate the process of the weights as
time passes.
\subsection{The weight of the first block as an O.U. type process}
In this section, we focus on the first block $\Pi_1(t)$, that is the
cluster at time $t$ which contains the root $1$ of $\T$.  The next
statement gathers its key properties, and in particular stresses the
connection with an Ornstein-Uhlenbeck type process.

\begin{theorem}\label{T2} 
{\rm (i)} For every $t\geq 0$, the following limit
$$\lim_{n\to\infty} n^{-\e^{-t}}\#\{j\leq n: j\in \Pi_1(t)\}=X_1(t)$$
exists in $(0,\infty)$ a.s. The variable $X_1(t)$ has the Mittag-Leffler
distribution with parameter $\e^{-t}$,
$$\P(X_1(t)\in \d x)/\d x=\frac{\e^t}{\pi} \sum_{k=0}^{\infty} \frac{(-1)^{k+1}}{k!} \Gamma(k\e^{-t}+1) x^{k-1}\sin(\pi k \e^{-t});$$
equivalently, its Mellin transform is given by
$$\E(X_1^q(t))= \frac{\Gamma(q+1)}{\Gamma(\e^{-t}q+1)}, \qquad q\geq 0.$$

\noindent{\rm (ii)} The process $(X_1(t): t\geq 0)$ is Markovian; its
semigroup $P_t(x,\cdot)$ is given by
$$P_t(x,\cdot)=\P(x^{\e^{-t}}X_1(t)\in\cdot).$$

\noindent{\rm (iii)} The process $Y(t)=\ln X_1(t)$, $t\geq 0$, is of
Ornstein-Uhlenbeck type. More precisely,
$$L(t)=Y(t)+\int_0^tY(s)\d s,\qquad t\geq 0,$$
is a spectrally negative L\'evy process with cumulant-generating function
$$\kappa(q)=\ln\E(\exp(qL(1))),\qquad q\geq 0,$$
given by
$$\kappa(q) = q\psi(q+1),$$ where $\psi$ denotes the digamma function, that
is the logarithmic derivative of the gamma function. 
\end{theorem}

 In the sequel, we shall refer to $X_1(t)$ as the weight of the first block
(or the root-cluster) at time $t$.  
 Before tackling the proof of Theorem~\ref{T2}, we make a couple of comments.

 Firstly, observe from (i) that $\lim_{t\to\infty}
 \E(X_1(t)^q)=\Gamma(q+1)$, so that as $t\to\infty $, $Y(t)$ converges in
 distribution to the logarithm of a standard exponential variable. On the
 other hand, it is well-known that the weak limit at $\infty$ of an
 Ornstein-Uhlenbeck type process is self-decomposable; cf. Section 17 in
 Sato~\cite{Sato}. So (iii) enables us to recover the fact that the
 log-exponential distribution is self-decomposable; see Shanbhag and
 Sreehari~\cite{SS}. 
 
 Secondly, note  that the L\'evy-Khintchin formula for $\kappa$ reads 
$$\kappa(q)= -\gamma q + \int_{-\infty}^{0}(\e^{qx}-1-qx)\frac{\e^{x}}{(1-\e^x)^2}\d x,$$
where $\gamma=0.57721\ldots$ is the Euler-Mascheroni constant. Indeed, this follows
readily from the classical identity for the digamma function
$$\psi(q+1)=-\gamma + \int_0^{1}\frac{1-x^q}{1-x} \d x.$$
In turn, this enables us to identify the L\'evy measure of $L$ as
$$\Lambda(\d x)={\e^{x}}{(1-\e^x)^{-2}}\d x,\qquad  x\in(-\infty, 0).$$
Since the jumps of $L$ and of $Y$ coincide, the L\'evy-It\=o decomposition
entails that the jump process of $Y= \ln X_1$ is a Poisson point process
with characteristic measure $\Lambda$.  In this direction, recall from
Proposition~\ref{P3} that the distribution of the asymptotic frequency of
the first block under the measure ${\bf r}$ of the splitting rates of $\Pi$
is $(1-y)^{-2} \d y$, $y\in(0,1)$, and observe that the image of the latter
by the map $y\mapsto \ln y$ is precisely $\Lambda$. This should of course
not come as a surprise.

We shall present two proofs of Theorem~\ref{T2}(i); the first relies on the
well-known connection between random recursive trees and Yule processes and
is based on arguments due to Pitman. Indeed, $\#\{j\leq n: j\in \Pi_1(t)\}$
can be interpreted in terms of the two-type population system considered in
Section 3.4 of~\cite{PiSF}, as the number of novel individuals at time $t$
when the birth rate of novel offspring per novel individual is given by
$\alpha=\e^{-t}$, and conditioned that there are $n$ individuals in total
in the population system at time $t$. Part (i) of the theorem then readily
follows from Proposition 3.14 in connection with Corollary 3.15 and Theorem
3.8 in~\cite{PiSF}. For the reader's convenience, let us nonetheless give a
self-contained proof which is specialized to our situation. We further
stress that variations of this argument will be used in the proofs of
Proposition~\ref{P5} and Corollary~\ref{C2}.
  
{\noindent{\bf First proof of Theorem~\ref{T2}(i):}\hskip10pt} Consider a
population model started from a single ancestor, in which each individual
gives birth a new child at rate one (in continuous time). If the ancestor
receives the label $1$ and the next individuals are labeled $2,3,\ldots$
according to the order of their birth times, then the genealogical tree of
the entire population is a version of $\T$. Further, if we write $Z(s)$ for
the number of individuals in the population at time $s$, then the process
$(Z(s): s\geq 0)$ is a Yule process, that is a pure birth process with
birth rate $n$ when the population has size $n$. Moreover, it is readily
seen that the Yule process $Z$ and the genealogical tree $\T$ are
independent.

It is well-known that 
$$\lim_{s\to\infty} \e^{-s}Z(s) = W\qquad \hbox{almost surely},$$
where $W$ has the standard exponential distribution. As a consequence, if
we write $\tau_n=\inf\{s\geq 0: Z(s)=n\}$ for the birth-time of the
individual with label $n$, then
\begin{equation}\label{E1}
\lim_{n\to\infty} n\e^{-\tau_n}=W \qquad \hbox{almost surely.}
\end{equation}

Now we incorporate destruction of edges to this population model by killing
merciless each new-born child with probability $1-p\in(0,1)$, independently
of the other children. The resulting population model is again a Yule
process, say $Z^{(p)}=(Z^{(p)}(s):s\geq 0)$, but now the rate of birth per
individual is $p$. Therefore, we have also
$$\lim_{s\to\infty} \e^{-ps}Z^{(p)}(s) = W^{(p)}\qquad \hbox{almost surely,}$$
where $W^{(p)}$ is another standard exponential variable. We stress that
$W^{(p)}$ is of course correlated to $W$ and not independent of $\T$,
in contrast to $W$.

In this framework, we identify for $p=\e^{-t}$
$$\#\{j\leq n: j\in \Pi_1(t)\}= Z^{(p)}(\tau_n)$$
and therefore
$$\lim_{n\to\infty} \e^{-p\tau_n} \#\{j\leq n: j\in
\Pi_1(t)\}=W^{(p)}\qquad \hbox{almost surely}.$$
Combining with~\eqref{E1}, we arrive at 
$$\lim_{n\to\infty} n^{-p} \#\{j\leq n: j\in
\Pi_1(t)\}=\frac{W^{(p)}}{W^p}\qquad \hbox{almost surely},$$
which proves the first part of (i).

Now recall that the left-hand side above only depends on the genealogical
tree $\T$ and the exponential random variables $\epsilon_i$ attached to its
edges. Therefore, it is independent of the Yule process $Z$ and {\it a
  fortiori} of $W$. Since both $W$ and $W^{(p)}$ are standard exponentials,
the second part of (i) now follows from the moments of exponential random
variables. \QED

The second proof of Theorem~\ref{T2}(i) relies on more advanced features on
the destruction of random recursive trees and Poisson-Dirichlet partitions.

{\noindent{\bf Second proof of Theorem~\ref{T2}(i):}\hskip10pt} It is known
from the work of Goldschmidt and Martin~\cite{GM} that the destruction of
$\T$ bears deep connections to the Bolthausen-Sznitman coalescent. In this
setting, the quantity
$$\#\{j\leq n: j\in \Pi_1(t)\}$$
can be viewed as the number of blocks at time $t$ in a Bolthausen-Sznitman
coalescent on $[n]=\{1,\ldots, n\}$ started from the partition into
singletons. On the other hand, it is known that the latter is a so-called
$(\e^{-t},0)$ partition; see Section 3.2 and Theorem 5.19 in
Pitman~\cite{PiSF}. Our claims now follow from Theorem 3.8
in~\cite{PiSF}. \QED

{\noindent{\bf Proof of Theorem~\ref{T2}(ii):}\hskip10pt} Let $\Pi'_1$ be
an independent copy of the process $\Pi_1$. Fix $s,t\geq 0$ and put
$B=\Pi_1(s)$, $C=\Pi'_1(t)$. Recall that $B(j)$ denotes the $j$-th smallest
element of $B$, and $B(C)$ stands for the block $\{B(j):j\in C\}$. By
Proposition~\ref{P1}, there is the equality in distribution
$\Pi_1(s+t)=B(C)$.  From (i) we deduce that
$$B(n)\sim \left(n/X_1(s)\right)^{\e^s}\qquad\hbox{almost surely as }n\to \infty,$$ and
similarly $C(n)\sim \left(n/X'_1(t)\right)^{\e^t}$ as $n\to\infty$, where
$X'_1(t)$ has the same law as $X_1(t)$ and is further independent of
$(X_1(r): r\geq 0)$. It follows that there are the identities
\pagebreak
\begin{align*}
X_1(s+t)&= \lim_{m\rightarrow\infty}m^{-\e^{-(s+t)}}\#\{j\leq m :j\in\Pi_1(s+t)\}\\
&=\lim_{n\rightarrow\infty}\left((B(C)(n))^{-\e^{-(s+t)}}n\right)\\
&=\lim_{n\rightarrow\infty}\left((B(C(n)))^{-\e^{-(s+t)}}n\right)\\
&=X_1^{\e^{-t}}(s)X'_1(t).
\end{align*}
Here, in the next to last equality we have used the fact that the $n$-th
smallest element of $B(C)$ is given by the $C(n)$-th smallest element of $B$,
and for the last equality we have plugged in the asymptotic
expressions for $B(n)$ and $C(n)$ that we found above. Our claim now follows easily.  \QED

We point out that, alternatively, the Markov property of $X_1$ can also be
derived from the interpretation of $\#\{j\leq n: j\in \Pi_1(t)\}$ as the
number of blocks at time $t$ in a Bolthausen-Sznitman coalescent on $[n]$;
see the second proof of Theorem~\ref{T2}(i) above.

{\noindent{\bf Proof of Theorem~\ref{T2}(iii):}\hskip10pt} We first observe
from (ii) that the process $Y$ is Markovian with semigroup $Q_t(y,\cdot)$
given by
$$Q_t(y, \cdot) = \P(\e^{-t}y + Y(t)\in \cdot).$$
Next, recall from the last remark made after Theorem~\ref{T2} that the
function $q \mapsto \kappa(q) = q\psi(q+1)$ is the cumulant-generating
function of a spectrally negative L\'evy process, say $L=(L(t): t\geq
0)$. Consider then the Ornstein-Uhlenbeck type process $U=(U(t): t\geq 0)$
that solves the stochastic differential equation
$$U(t) = L(t) - \int_0^tU(s)\d s,$$
that is, equivalently, $U(t) = \e^{-t}\int_0^t\e^s \d L(s)$. Then $U$ is
also Markovian with semigroup $R_t(u,\cdot)$ given by
$$R_t(u, \cdot) = \P(\e^{-t}u + U(t)\in \cdot).$$
So to check that the processes $Y$ and $U$ have the same law, it suffices
to verify that they have the same one-dimensional distribution.

The calculations of Section 17 in Sato~\cite{Sato} (see Equation (17.4) and
Lemma 17.1 there) show that for every $q\geq 0$,
$$\E\left(\exp (qU(t))\right) = \exp\left(\int_0^t\kappa(\e^{-s}q )\d s\right).$$
Now observe that 
 $$\int_0^t\kappa(\e^{-s}q )\d s= \int_0^t \e^{-s}q \psi(\e^{-s}q +1)\d s = \ln \Gamma(q+1) - \ln \Gamma(\e^{-t}q+1),$$
 and hence
$$\E\left(\exp (qU(t))\right)= \frac{\Gamma(q+1)}{\Gamma(\e^{-t}q+1)}= \E\left(\exp (qY(t))\right),$$
where the second identity stems from Theorem~\ref{T2}(i). \QED

\noindent{\bf Remark.}
It would be interesting to understand whether with probability $1$, the
block $\Pi_1(t)$ has weights in the sense of Theorem~\ref{T2}(i)
simultaneously for all $t\geq 0$ (note that the asymptotic frequencies are
equal to zero on $(0,\infty)$ a.s.), and whether $t\mapsto X_1(t)$ is
c\`adl\`ag.  In this direction, we recall that each block of a standard
homogeneous fragmentation process possesses asymptotic frequencies
simultaneously for all $t\geq 0$ a.s., see Proposition 3.6 in~\cite{RFCP}.
Moreover, if $B$ denotes a block of such a process, then the process
$t\mapsto|B(t)|$ is c\`adl\`ag.  Here it should be observed that the first
block $\Pi_1$ is the only block which is decreasing, in the sense that
$\Pi_1(t')$ is contained in $\Pi_1(t)$ a.s. whenever $t'\geq t$. This
however does not imply monotonicity of $X_1(t)$ in $t$, which is crucial
ingredient for the proof of the mentioned properties in the case of a
homogeneous fragmentation.

\subsection{Fragmentation of weights as a branching O.U. process}
We next turn our interest to the other blocks of the partition $\Pi(t)$; we
shall see that they also have a weight, in the same sense as for the first
block.  In this direction, it is convenient to write first $\T_i$ for the
subtree of $\T$ rooted at $i\geq 1$; in particular $\T_1=\T$. Then for
$t\geq 0$, we write $T_i(t)$ the subtree of $\T_i$ consisting of vertices
$j\in\T_i$ which are still connected to $i$ after the edges $e_k$ with
$\epsilon_k\leq t$ have been removed.  Note that for $i\geq 2$, $T_i(t)$ is
a cluster at time $t$ if and only if $\epsilon_i\leq t$, an event which has
always probability $1-\e^{-t}$ and is further independent of $T_i(t)$.  On
that event, the vertex set of $T_i(t)$ is a block of the partition
$\Pi(t)$, and all the blocks of $\Pi(t)$ arise in this form.

\begin{lemma}\label{L1}
 For every $t\geq 0$ and $i\in\N$, the following limit
$$\lim_{n\to\infty} n^{-\e^{-t}}\#\{j\leq n: j\in T_i(t)\}=\rho_i(t)$$
exists in $(0,\infty)$ a.s. Moreover, the process 
$$\rho_i=(\rho_i(t):t\geq 0)$$
has the same law as 
$$(\beta_i^{\e^{-t}}X_1(t): t\geq 0),$$
where $\beta_i$ denotes a beta variable with parameter $(1,i-1)$ and is further independent of $X_1(t)$.
In particular, the positive moments of $\rho_i(t)$ are given by 
$$\E(\rho_i^q(t))= \frac{\Gamma(q+1)\Gamma(i)}{\Gamma(\e^{-t}q+i)}, \qquad q\geq 0.$$
\end{lemma}
\proof The recursive construction of $\T$ and $\T_i$ has the same dynamics
as a P\'olya urn, and basic properties of the latter entail that the
proportion $\beta_i$ of vertices in $\T_i$ has the beta distribution with
parameter $(1,i-1)$. Further, enumerating the vertices of $\T_i$ turns the
latter into a random recursive tree. Our claim then follows readily from
Theorem~\ref{T2}. \QED

Lemma~\ref{L1} entails that for every $i\in\N$, the $i$-th block $\Pi_i(t)$
of $\Pi(t)$ has a weight in the sense of~\eqref{eqmass}, a.s.  We write
$X_i(t)$ for the latter and set $\X(t)=(X_1(t), X_2(t), \ldots)$. We now
investigate the process $\X=(\X(t): t\geq 0)$.

Firstly, using the functional equation of the gamma function, the
integral representation of the beta function and the expression for the
moments of $X_1(t)$ from Theorem~\ref{T2}(i), an easy calculation shows
\begin{equation}
\label{eq-moments}
\E\left(\sum_{i=1}^{\infty}X_i^q(t)\right) = \E\left(X_1^q(t)\right) +
\E\left(\sum_{i=2}^{\infty}\1_{\{\epsilon_i\leq t\}}\rho_i^q(t)\right)= 
\frac{(q-1)}{(\e^{-t}q-1)}\frac{\Gamma(q)}{\Gamma(\e^{-t}q)},
\end{equation}
provided $q>\e^t$. In particular, for every $t\geq 0$, the $X_i(t)$ can be
sorted in the decreasing order. We write $\X^{\downarrow}(t)$ for the
sequence obtained from $\X(t)$ by ranking the weights $X_i(t)$
decreasingly, where as usual elements are repeated according to their
multiplicity. For $q>0$, let
$$\ell^{q\downarrow}=\left\{{\bf x} = (x_1,\dots):x_1\geq x_2\geq\dots\geq
  0,\hbox{ and }\sum_{i=1}^\infty x_i^q<\infty\right\},
$$
endowed with the $\ell^q$-distance. Similarly, denote by
$\ell^{\infty\downarrow}$ the space of ordered sequences of
positive reals, endowed with the $\ell^\infty$-distance. 
For the process $\X^\downarrow$, we obtain the following characterization.
\begin{corollary}\label{C1}
  Let $T\in(0,\infty]$, and set $q=\e^{T}$ (with the convention
  $\e^\infty=\infty$).  Then the process $\X^\downarrow=(\X^\downarrow(t) :
  t< T)$ takes its values in $\ell^{q\downarrow}$ and is Markovian. More
  specifically, its semigroup can be described as follows. For $s,t \geq
  0$ with $s+t<T$, the law of $\X^\downarrow(s+t)$ conditioned on
  $\X^\downarrow(s)= (x_1,\dots)$ is given by the distribution of the
  decreasing rearrangement of the sequence $(x_i^{\e^{-t}}x_j^{(i)} :
  i,j\in\N)$, where $((x_1^{(i)},\dots):i\in\N)$ is a
  sequence of independent random elements in $\ell^{q\downarrow}$, each of
  them distributed as $\X^{\downarrow}(t)$.
\end{corollary}
\proof The fact that $\X^\downarrow(t)\in\ell^{q\downarrow}$ for $t<T$
follows from~\eqref{eq-moments}. The specific form of the semigroup can be
deduced from Theorem~\ref{T1} and a variation of the arguments leading to
Theorem~\ref{T2}(ii).  \QED

We now draw our attention to the logarithms of the weights. In this
direction, recall that for a homogeneous fragmentation, the process of the
asymptotic frequencies of the blocks bears close connections with branching
random walks. More precisely, the random point process with atoms at the
logarithm of the asymptotic frequencies and observed, say at integer times,
is a branching random walk; see~\cite{BR} and references therein. This
means that at each step, each atom, say $y$, is replaced by a random cloud
of atoms located at $y+z$ for $z\in {\mathcal Z}$, independently of the
other atoms, and where the random point process ${\mathcal Z}$ has a fixed
distribution which does not depend on $y$ nor on the step. In the same
vein, we also point out that recently, a natural extension of homogeneous
fragmentations, called compensated fragmentations, has been constructed
in~\cite{Be3}, and bears a similar connection with branching L\'evy
processes. Note that compensated fragmentations are
$\ell^{2\downarrow}$-valued processes, in contrast to $\X^\downarrow$.

Our observations incite us to introduce
$$\pY_t = \sum_{i=1}^{\infty} \delta_{\ln X_i(t)},\qquad t\geq 0.$$

\begin{theorem}\label{T3} The process with values in the space of point
  measure $\pY=(\pY_t : t\geq 0)$ is an Ornstein-Uhlenbeck branching
  process, in the sense that it is Markovian and its transition
  probabilities can be described as follows:

  \noindent
For every $s,t\geq 0$, the conditional law of $\pY_{s+t}$ given 
$\pY_s=\sum_{i=1}^{\infty} \delta_{y_i}$ is given by the distribution of
$$\sum_{i=1}^{\infty} \sum_{j=1}^{\infty} \delta_{\e^{-t}y_i+\zeta^{(i)}_{j}}$$
where the point measures 
$$\pZ^{(i)} = \sum_{j=1}^{\infty} \delta_{\zeta^{(i)}_{j}}$$
are independent and each has the same law as $\pY_t$. 

Furthermore, the mean intensity of $\pY_t$ is determined by 
$$\E\left( \int \e^{qy}\pY_t(\d y)\right) = \frac{(q-1)}{(\e^{-t}q-1)}\frac{\Gamma(q)}{\Gamma(\e^{-t}q)},\qquad q>\e^t.$$
\end{theorem}

\proof The first claim follows from Theorem~\ref{T1} and an adaption of the
argument for Theorem~\ref{T2}(iii). For the second, we note that 
$$
\E\left( \int \e^{qy}\pY_t(\d y)\right)=\E\left(\sum_{i=1}^{\infty}X_i^q(t)\right),
$$
and the right hand side has already been evaluated in~\eqref{eq-moments}.
\QED

Next we give a description of the finite dimensional laws of $\X(t) =
(X_1(t),X_2(t),\ldots)$ for $t>0$ fixed. In this direction, it is
convenient to define two families of probability distributions.

The first family is indexed by $j\in \N$ and $t>0$, and  is defined as 
$$
\mu_{j,t}(k) = {k-2\choose
  k-j-1}(\e^{-t})^{k-j-1}(1-\e^{-t})^{j},\qquad k\geq j+1.
$$
Note that the shifted distribution $\tilde{\mu}_{j,t}(k) = \mu_{j,t}(k+1)$,
$k\geq j$, is sometimes called the negative binomial distribution with
parameters $j$ and $1-\e^{-t}$, that is the law of the number of
independent trials for $j$ successes when the success probability is given
by $1-\e^{-t}$.

The second family is indexed by $j\in\N$ and $k\geq j$, and can be
described as follows. We denote by $\theta_{j,k}$ the probability measure
on the discrete simplex $\Delta_{k,j}=\{(k_1, \ldots, k_j)\in\N^j:
k_1+\cdots + k_j=k\}$, such that
$\theta_{j,k}\left(k_1,\ldots,k_{j}\right)$ is the probability that on a
random recursive tree of size $k$ (that is on a random tree distributed as
$\T_{\mid[k]}$), after $j-1$ of its edges chosen uniformly at random have
been removed, the sequence of the sizes of the $j$ subtrees, ordered
according to the label of their root vertex, is given by
$(k_1,\ldots,k_{j})$.

\noindent{\bf Remark.}
 The distribution $\theta_{j,k}$ is
equal to $\delta_k$ for $j=1$. For
$j=2$, Meir and Moon~\cite{MM} found the expression
$$
\theta_{2,k}(k_1,k_2) = \frac{k}{k_2(k_2+1)(k-1)},\qquad k_1,k_2\in\N
\hbox{ with }k_1+k_2=k,
$$
with $\theta_{2,k}(k_1,k_2)=0$ for all other pairs
$(k_1,k_2)$. Generalizing the proof of this formula given in~\cite{MM} to
higher $j$, we find $(k\geq j\geq 3$ and $k_1+\dots+k_j=k$)
\begin{align*}
\lefteqn{\theta_{j,k}(k_1,k_2,\ldots,k_j)
=\frac{(k_1-1)!(k_2-1)!\cdots(k_j-1)!}{(k-1)!(k-1)\cdots(k-(j-1))}
\sum_{\ell_j=j-1}^{k-k_j}{k-\ell_j\choose
  k_j}}
\\
&\times
\sum_{\ell_{j-1}=j-2}^{(k-k_j-k_{j-1})\wedge(\ell_j-1)}{k-k_j-\ell_{j-1}\choose
  k_{j-1}}\times\dots\times
\sum_{\ell_2=1}^{(k-\sum_{i=2}^jk_i)\wedge(\ell_{3}-1)}{k-\sum_{i=3}^{j}k_i-\ell_2\choose k_2}.
\end{align*}

\begin{proposition}
\label{P5}
  Let $j\in\N$, $q_1,\ldots,q_{j+1} \geq 0$, and set $k_{j+1}=1$. The
  Mellin transform of the vector $(X_1(t),\dots X_{j+1}(t))$ for fixed
  $t>0$ is given by
$$
\E\left(X^{q_1}_1(t)\cdots X^{q_{j+1}}_{j+1}(t)\right) =
\sum_{k=j+1}^\infty\mu_{j,t}(k)\sum_{k_1,\ldots,k_{j}\geq 1,\atop k_1+\dots +k_{j}=k-1}
\theta_{j,k-1}\left(k_1,\ldots,k_{j}\right)\frac{\Gamma(k)}{\Gamma(q\e^{-t}+k)}\prod_{i=1}^{j+1}\frac{\Gamma(q_i+k_i)}{\Gamma(k_i)}.
$$
\end{proposition}
\noindent{\bf Remark.}
By plugging in the definition of $\mu_{j,t}(k)$, one checks that the
right hand side is finite.

\proof Fix $t>0$, and set $p=\e^{-t}$. For ease of notation, we write
$\Pi_i$ and $X_i$ instead of $\Pi_i(t)$ and $X_i(t)$.  Furthermore, fix an
integer $j\in\N$ and numbers $k_1,\ldots,k_{j}\in\N$. For convenience, set
$k=k_1+\dots+k_j+k_{j+1}$, with $k_{j+1}=1$. We first work conditionally on
the event
$$
A_{k_1,\dots,k_j}=\big\{\min \Pi_{j+1} = k,\, \#(\Pi_1\cap[k]) =
k_1,\dots,\#(\Pi_j\cap[k]) = k_j\big\}.
$$
We shall adapt the first proof of
Theorem~\ref{T2}(i). Here, we consider a multi-type Yule process starting
from $k$ individuals in total such that $k_i$ of them are of type $i$, for
each $i=1,\ldots,j+1$. The individuals reproduce independently of each
other at unit rate, and each child individual adopts the type of its
parent. Then, if $Z(s)$ stands for the total number of individuals at time
$s$, we have that $\lim_{s\rightarrow\infty}\e^{-s}Z(s)=\gamma(k)$ almost
surely, where $\gamma(k)$ is distributed as the sum of $k$ standard
exponentials, i.e. follows the gamma law with parameters $(k,1)$. Now
assume again that each new-born child is killed with probability
$1-p\in(0,1)$, independently of each other. Writing $Z^{(i,p)}(s)$ for the
size of the population of type $i$ at time $s$ (with killing), we obtain 
$$
\lim_{s\rightarrow\infty}\e^{-ps}Z^{(i,p)}(s)=\gamma^{(i,p)}(k_i),\qquad i=1,\dots,j+1,
$$
where the $\gamma^{(i,p)}(k_i)$ are independent
gamma$(k_i,1)$ random variables (they are however clearly correlated to
the asymptotic total population size $\gamma(k)$). From the arguments
given in the first proof of Theorem~\ref{T2}(i) it should be plain that
conditionally on the event $A_{k_1,\dots,k_j}$, we have for the weights
$X_i$ the representation
$$
X_i = \frac{\gamma^{(i,p)}(k_i)}{\gamma^p(k)},\qquad i=1,\dots,j+1,
$$
and the $X_i$ are independent of $\gamma(k)$.  Now let
$q_1,\ldots,q_{j+1}\geq 0$ and put $q=q_1+\dots+q_{j+1}$. Using the
expression for the $X_i$ and independence, we calculate
\begin{equation}\label{eq56}
\E\left(\gamma(k)^{qp}\right)\E\left(X^{q_1}_1\cdots
  X^{q_{j+1}}_{j+1}\mid A_{k_1,\dots,k_j}\right) =\prod_{i=1}^{j+1}\frac{\Gamma(q_i+k_i)}{\Gamma(k_i)}.
\end{equation}
Therefore, again with $k_{j+1}=1$,
$$
\E\left(X^{q_1}_1\cdots
  X^{q_{j+1}}_{j+1}\right) = \sum_{k=j+1}^\infty\frac{\Gamma(k)}{\Gamma(qp+k)}\sum_{k_1,\ldots,k_{j}\geq 1,\atop
  k_1+\dots+ k_{j}=k-1}\prod_{i=1}^{j+1}\frac{\Gamma(q_i+k_i)}{\Gamma(k_i)}\P\left(A_{k_1,\dots,k_j}\right).
$$
With $k=k_1+\dots+k_{j+1}$ as above, we express the probability of $A_{k_1,\dots,k_j}$ as
$$
\P\left(A_{k_1,\dots,k_j}\right)= \P\left(\#(\Pi_1\cap[k]) =
k_1,\dots,\#(\Pi_j\cap[k]) = k_j\mid \min \Pi_{j+1} = k\right)\P\left(\min \Pi_{j+1} = k\right).
$$
By induction on $j$ we easily deduce that $\min \Pi_{j+1}-1$ is distributed as the
sum of $j$ independent geometric random variables with success
probability $1-p$, i.e. $\min \Pi_{j+1}-1$ counts the number of trials for
$j$ successes, so that $\P\left(\min \Pi_{j+1} =
  k\right)=\mu_{t,j}(k)$. Moreover, it follows from the very definition of
the blocks $\Pi_i$ and the fact that the exponentials attached to the edges
of $\T_{\mid[k]}$ are i.i.d. that
$$\P\left(\#(\Pi_1\cap[k]) =
k_1,\dots,\#(\Pi_j\cap[k]) = k_j\mid \min \Pi_{j+1} = k\right) =
\theta_{j,k-1}(k_1,\dots,k_j).$$
This proves the proposition.
\QED

We finally look closer at the joint moments of
$X_1(t)$ and $X_2(t)$ when $t$ tends to zero. We observe
$\theta_{1,k}=\delta_k$ and $\mu_{1,t}(k)=(\e^{-t})^{k-2}(1-\e^{-t})$, so
that
$$
\E\left(X^{q_1}_1(t)X^{q_2}_2(t)\right)=(1-\e^{-t})\Gamma(q_2+1)\sum_{k=2}^\infty(k-1)(\e^{-t})^{k-2}\frac{\Gamma(q_1+k-1)}{\Gamma((q_1+q_2)\e^{-t}+k)}.
$$
Now assume $q_2>1$. From the last display we get 
\begin{align*}
  \lim_{t\rightarrow 0+}\frac{1}{t}\E\left(X_1^{q_1}(t)X_2^{q_2}(t)\right)
  &=\sum_{k=1}^\infty
  k\frac{\Gamma(q_1+k)\Gamma(q_2+1)}{\Gamma(q_1+q_2+k+1)}
  \\
  &=\sum_{k=1}^\infty k\int_0^1(1-x)^{q_1+k-1}x^{q_2}\dt x= \int_0^1
  (1-x)^{q_1}x^{q_2-2}\dt x,
\end{align*}
which one could have already guessed from Proposition~\ref{P3}.

\subsection{Asymptotic behaviors}
We shall finally present some asymptotic properties of the process $\X$ of
the weights.  To start with, we consider the large time behavior.

\begin{corollary} \label{C2} 
As $t\to \infty$, there is the weak convergence
$$(X_i(t): i\in \N) \Longrightarrow (W_i: i\in\N),$$
where on the right-hand side, the $W_i$ are i.i.d. standard exponential
variables.
\end{corollary}

\noindent{\bf Remark.} This result is a little bit surprising, as obviously
$\Pi(\infty)$ is the partition into singletons. That is $\Pi_i(\infty)$ is
reduced to $\{i\}$ and hence has weight $1$ if we apply \eqref{eqmass} for
$t=\infty$. In other words, the limits $n\rightarrow\infty$ and
$t\rightarrow\infty$ may not be interchanged.

\proof Fix $j\in\N$ arbitrarily and consider the event $A(t)=\{\min
\Pi_{j+1}(t)=j+1\}$. Recall from the proof of Proposition \ref{P5} that
$\P(A(t))=\mu_{t,j}(j+1)$ and note that this quantity converges to $1$ as
$t\to \infty$. Further, on the event $A(t)$, we have also $\Pi_i(t)\cap
[j+1]=\{i\}$ for all $1\leq i \leq j+1$, that is $A(t)=A_{1, \ldots, 1}$,
again in the notation of the proof of Proposition \ref{P5}.

Take $q_1, \ldots, q_{j+1}\geq 0$. Applying \eqref{eq56}, we get
$$\lim_{t\to\infty}\E\left(X^{q_1}_1(t)\cdots X^{q_{j+1}}_{j+1}(t)\mid A(t)\right) = \prod_{i=1}^{j+1}\Gamma(q_i+1),
$$
which entails that $(X_1(t), \ldots, X_{j+1}(t))$ converge in distribution
as $t\to \infty$ towards a sequence of $j+1$ independent exponentially
distributed random variables. \QED

We next consider for $t>0$ fixed the behavior of $X_n(t)$ as $n\to
\infty$.

\begin{corollary}\label{C3} Let $t>0$. As $n\to \infty$, there is the weak
  convergence
$$n^{\e^{-t}} X_n(t)  \Longrightarrow V^{\e^{-t}}X_1(t),$$
where $V$ denotes an exponential random variable of parameter
$(1-\e^{-t})^{-1}$ which is independent of $X_1(t)$.
\end{corollary}

\proof In the notation of Lemma~\ref{L1}, we have
$X_n(t)=\rho_{i(n,t)}(t)$, where (again with the convention $\epsilon_1=0$)
$$i(n,t)=\min\left\{j\geq 1: \sum_{i=1}^j\1_{\{\epsilon_i\leq t\}}=n\right\}.$$
From Lemma~\ref{L1} we know that $X_n(t)=\beta^{\e^{-t}}_{i(n,t)}X_1(t)$ in
distribution, where $i(n,t)$ and $\beta_{i(n,t)}$ are both independent of
$X_1(t)$. By the law of large numbers, $i(n,t)\sim n(1-\e^{-t})^{-1}$
almost surely. Writing
$$\beta_{i(n,t)}=\left(i(n,t)\beta_{i(n,t)}\right)i(n,t)^{-1}$$
and using the fact that $k\beta_k$ converges in distribution to a standard
exponential random variable as $k\rightarrow\infty$, the claim follows.
\QED

We finally look at the case $t\rightarrow 0+$. From Theorem~\ref{T2}(i) it
follows that
\begin{equation}
\label{E2}
\lim_{t\rightarrow 0+}X_1(t) = 1\qquad\hbox{in probability.}
\end{equation}
In fact, $X_1(t)$ tends to be the largest element of the sequence $\X$ when
$t$ approaches zero.
\begin{corollary}\label{C} 
$$
\lim_{t\rightarrow 0+}\P\left(X_1(t)>X_i(t)\textup{ for all }i\geq
  2\right) = 1.
$$
\end{corollary}
\proof
We consider the complementary event. We have
$$
\P\left(X_i(t)\geq X_1(t)\textup{ for some }i\geq
  2\right) \leq \P\left(X_i(t)\geq 1/2 \textup{ for some }i\geq
  2\right) + \P\left(X_1(t)\leq 1/2\right),
$$
and the second probability on the right hand side converges to zero as
$t\rightarrow 0+$ by~\eqref{E2}. For the first probability, we have by
Lemma~\ref{L1}, with $\beta_i$ denoting a beta$(1,i-1)$ random variable,
$$
\P\left(X_i(t)\geq 1/2 \textup{ for some }i\geq
  2\right) \leq
(1-\e^{-t})\sum_{i=2}^\infty\P\left(\beta_i X_1^{\e^{t}}(t)\geq (1/2)^{\e^{t}}\right).
$$
Using independence and the expression for the moments of $X_1(t)$ from
Theorem~\ref{T1}, we obtain for $t\geq 0$ such that $\e^t\leq 2$,
$\E(\beta_i X_1^{\e^{t}}(t))\leq  2/i$ and Var$(\beta_i
X_1^{\e^{t}}(t))\leq 24/i^2.$ Therefore, for such $t$ and $i\geq 10$, by
Chebycheff's inequality, with $C=10^4$,
$$
\P\left(\big|\beta_i X_1^{\e^{t}}(t)-\E(\beta_i X_1^{\e^{t}}(t))\big|\geq
  (1/2)^{\e^{t}}-2/i\right) \leq \frac{C}{i^{2}}.
$$
This shows
$$
\P\left(X_i(t)\geq 1/2 \textup{ for some }i\geq
  2\right) \leq
(1-\e^{-t})\left(8+C\sum_{i=10}^\infty
  \frac{1}{i^2}\right)=O(t)\qquad\hbox{as }t\rightarrow 0+.
$$
\QED

\noindent {\bf Concluding remark.} 
As it should be plain from the introduction, the set of vertices
$$C_i^{(n)}(t)=\{j\leq n : j\in\Pi_i(t)\}$$
form the percolation clusters of a Bernoulli bond percolation with
parameter $p=\e^{-t}$ on a random recursive tree of size $n$ on the vertex
set $\{1,\dots,n\}$. Our results of Section~\ref{SW} can therefore be
understood as results on the asymptotic sizes of these clusters when $n$
tends to infinity.

Cluster sizes of random recursive trees were already studied in~\cite{Be2}
when the percolation parameter $p$ satisfies $p=p(n)=1-s/\ln n +o(1/\ln n)$ for
$s>0$ fixed. The analysis was extended in~\cite{Ba} to all regimes
$p(n)\rightarrow 1$. It shows that in these regimes, the root
cluster containing $1$ has always the asymptotic size $\sim n^{p(n)}$, while the
next largest cluster sizes, normalized by a factor
$(1-p(n))^{-1}n^{-p(n)}$, are in the limit given by the (ranked) atoms of a
Poisson random measure on $(0,\infty)$ with intensity $a^{-2}\dt a$. 

The regime of constant parameter $p=\e^{-t}$ considered here deserves the
name ``critical'', since in this regime, the root cluster and the next
largest clusters have the same order of magnitude, namely $n^{p}$. This is
already apparent from Lemma~\ref{L1}. More precisely, Corollary~\ref{C2}
readily shows that in fact $$\lim_{t\uparrow
  \infty}\liminf_{n\rightarrow\infty}\P\left(\textup{there exists }i\geq
  2\textup{ such that }\#C_i^{(n)}(t)> \#C_1^{(n)}(t)\right)=1.$$

\noindent {\bf Acknowledgment.}
The first author would like to thank Gr\'egory Miermont for stimulating
discussions.

\end{document}